\documentclass[12pt]{article}
\usepackage[cp1251]{inputenc}
\usepackage{amsmath,amssymb, amscd, amsthm}
\usepackage[english,russian]{babel}
\usepackage{graphicx}
\newcommand{\goe}{\geqslant}
\newcommand{\loe}{\leqslant}

\theoremstyle{plain}
\newtheorem{thm}{Теорема}
\newtheorem*{prop}{Предложение}
\newtheorem*{conja}{Гипотеза Калаи}
\newtheorem*{conjb}{Гипотеза B}
\newtheorem*{conj}{Гипотеза}

\theoremstyle{definition}

\newcommand{\Z}{\mathbb{Z}}

\DeclareMathOperator{\conv}{conv}
\DeclareMathOperator{\aff}{aff}
\DeclareMathOperator{\lin}{lin}

\begin{document}

\sloppy 

\title{\def\thefootnote{\arabic{footnote}} Совершенные призмоиды и решетчатые многогранники Делоне}
\author{\def\thefootnote{\arabic{footnote}} Марина Козачок\footnotemark[1], Александр Магазинов\footnotemark[2]}
\date{}

{\def\thefootnote{\arabic{footnote}} 
\footnotetext[1]{Математический институт им. Стеклова РАН;  marina.kozachok@gmail.com}
\footnotetext[2]{Математический институт им. Стеклова РАН;  magazinov-al@yandex.ru
}
\maketitle

\begin{abstract}
В данной работе доказано, что любой совершенный призмоид аффинно эквивалентен некоторому $0/1$-многограннику, а также некоторому решетчатому многограннику Делоне.
\end{abstract}
\section{Введение}
В 1989 году Г.Калаи сформулировал $3^d$-гипотезу:
\begin{conja} У любого выпуклого центрально-симметричного $d$-мерного многогранника общее число граней всех размерностей не меньше  $3^d$. \end{conja}В случае  $d=3$ доказательство гипотезы легко следует из формулы Эйлера, а равенство достигается только на кубе и дуальном ему октаэдре.  Для $d=4$ гипотеза была доказана в 2007 г. Г.Циглером, Р.Саньялом и А.Вернером в \cite{SWZ09}.  Верна ли $3^d$-гипотеза при $d\goe 5$ неизвестно. 

Нетрудно проверить, что у многогранников Ханнера (см. \cite{Ha56}) ровно $3^d$ граней. Например, для любой размерности многогранниками Ханнера являются куб и дуальный ему кроссполитоп. Г.Калаи предположил, что многогранники Ханнера исчерпывают класс многогранников, на которых достигается равенство в $3^d$-гипотезе. Однако этот вопрос остаётся открытым.

В пространстве размерности $d\loe 4$ доказана более сильная гипотеза (гипотеза В):
\begin{conjb}
Для любого центрально-симметричного многогранника $P$  существует многогранник Ханнера  $H$ такой, что $f_k(P)\goe f_k(H)$ для любого $k\,(0\loe k\loe d)$, где $f_k$ — число $k$-мерных граней соответсвующего многогранника.
\end{conjb} 
Однако, при  $d\goe 5$ гипотеза B  не верна (см.  \cite{SWZ09}).

В связи с $3^d$-гипотезой интерес представляет доказательство следующей ослабленной версии Гипотезы B:
\begin{conj}
Для любого центрально-симметричного многогранника $P$  существует центрально-симметричный совершенный призмоид $Q$ такой, что $f_k(P)\goe f_k(Q)$ для любого $k\,(0\loe k\loe d)$, где $f_k$ — число $k$-мерных граней соответсвующего многогранника.
\end{conj} 
Напомним, что центрально-симметричный совершенный призмоид — это центрально-симметричный многогранник, который является выпуклой оболочкой любой пары антиподальных граней. (Определение совершенного призмоида в общем случае дано ниже).  В  \cite{MK12} доказано, что любой многогранник Ханнера является совершенным призмоидом, но при $d\goe 5$ обратное уже не верно.

Основной результат данной работы — доказательство того, что любой совершенный призмоид аффинно эквивалентен некоторому решётчатому многограннику Делоне.  В \cite{Vo08}  доказано, что число аффинных типов решётчатых многогранников Делоне данной размерности конечно. Проведена классификация многогранников Делоне  для размерностей $d\loe 6$ (см. \cite{Du04}, \cite{Ba99}), что позволяет классифицировать и все совершенные призмоиды в этих размерностях.

\section{Определения}
Будем использовать следующие обозначения. Через $\conv X$ обозначим выпуклую оболочку множества точек $X$ в пространстве $\mathbb R^d$, через
$\aff X$ --- {\it аффинную оболочку}, т.е. множество всех линейных комбинаций 
\begin{multline*}
a_1x_1 + a_2x_2 + \ldots + a_nx_n, \quad \text{где\quad} 
a_1 + a_2 + \ldots + a_n = 1 \quad \text{и} \quad x_i \in X \quad (1\leq i \leq n).
\end{multline*}
Наконец, через $\lin X$ будем обозначать линейное пространство, ассоциированное с аффинной оболочкой $\aff X$, т.е. множество всех линейных комбинаций 
\begin{multline*}
a_1x_1 + a_2x_2 + \ldots + a_nx_n, \quad \text{где \quad}  
a_1 + a_2 + \ldots + a_n = 0 \quad \text{и} \quad x_i \in X \quad (1\leq i \leq n).
\end{multline*}

В данной работе рассматриваются многогранники в $d$-мерном евклидовом пространстве $\mathbb R^d$. Выпуклый $d$-мерный многогранник $P$ называется 
{\it призмоидом}, если некоторая пара параллельных $(d-1)$-мерных плоскостей содержит все его вершины. Примерами призмоидов являются, в частности, пирамида, 
призма и антипризма. 

Иначе говоря, призмоид $P$ есть выпуклая оболочка $\conv(F \cup F')$ двух многогранников $F$ и $F'$, каждый из которых лежит в одной из двух параллельных гиперплоскостей. При этом $\dim (\lin F + \lin F') = \dim P - 1$, где $\lin F + \lin F'$ понимается как сумма линейных подпространств в объемлющем пространстве.
Если $\dim P = d$ и $P \subset \mathbb R^d$, то две плоскости, параллельные $\lin F + \lin F'$ и проходящие через $F$ и $F'$, являются опорными для $P$.
Поэтому многогранники $F$ и $F'$ --- грани многогранника $P$. При этом каждая вершина призмоида $P$ принадлежит либо грани $F$, либо грани $F'$. 

Пусть $F$ --- гипергрань многогранника $P$, и существует такая грань $F'$, что  $\lin F'\subseteq \lin F$ и $P = \conv(F\cup F')$ (в дальнейшем будем говорить, что такие грани {\it параллельны}). Тогда будем называть $P$ 
{\it призмоидом над гранью $F$}.
 
Многогранник $P$ называется {\it совершенным призмоидом}, если он является призмоидом над любой своей гипергранью, т.е. для любой его гиперграни $F$ 
выполняется $P = \conv(F\cup F')$, где $F'$ --- параллельная к $F$ грань,  возможно, меньшей размерности.

Многогранник называется $0/1$-{\it многогранником}, если он является выпуклой оболочкой некоторого подмножества множества вершин единичного куба.

Пусть в евклидовом пространстве $\mathbb R^d$ задана некоторая решетка $\Lambda$. Пусть многогранник $D$ вписан в шар $B$,  причём внутренность шара $B$ не содержит ни одной точки решетки $\Lambda$ ($\partial B$ --- {\it пустая сфера} по Делоне), а также выполнено условие 
$$D = \conv( \Lambda \cap \partial B).$$
Тогда $D$ называется решетчатым многогранником Делоне (см. ~\cite{De37}). 
\section{Результаты}
\begin{prop}\label{prop:1}
Пусть $P$ --- совершенный призмоид c $N$ гипергранями. Тогда $P$ можно задать системой линейных неравенств
$$b_i\goe(a_i,x)\goe c_i, \qquad i = 1, 2, \ldots, N,$$
где  гиперплоскость, содержащая $i$-ю гипергрань, задается уравнением $(a_i,x) =  b_i$, и для любого $i = 1, 2, \ldots, N$ и любой вершины $v\in P$ верно либо $(a_i, v) =  b_i$, либо $(a_i, v) =  c_i$.
\end{prop}

\begin{proof}
Пусть гиперплоскость, содержащая гипергрань $F_i$, задается уравнением $(a_i, x) =  b_i$, при этом можно считать (за счет выбора знака $a_i$), что $P$ лежит в полупространстве
$(a_i, x) \leq b_i$. По определению совершенного призмоида, $P$ = $\conv (F_i \cup F'_i)$,
причем $\lin F'_i \subseteq \lin F_i$. Это значит, что грань $F'_i$ лежит в гиперплоскости $(a_i, x) =  c_i$, и $c_i < b_i$. 

Поскольку любая вершина $v$ призмоида $P$ принадлежит либо гиперграни $F$, либо грани $F'$, верно либо $(a_i, v) =  b_i$, либо $(a_i, v) =  c_i$.

Так как выпуклый многогранник определяется своими гипергранями, то система неравенств  $(a_i, x) \leq  b_i$ уже определяет многогранник $P$. С другой стороны, любая вершина $v\in Vert(P)$ (где $Vert(P)$ — множество вершин многоранника $P$), удовлетворяет
неравенству $(a_i, v) \geq  c_i$, поскольку $(a_i, v)$ принимает только значения $b_i$ и $c_i$. Поэтому добавление неравенств $(a_i, x) \geq  c_i$
не изменяет многогранник.

\end{proof}

\begin{thm}\label{01}
Для любого $d$-мерного совершенного призмоида существует аффинно эквивалентный ему $0/1$-многогранник, который получен из $d$-мерного единичного куба.
\end{thm}
\begin{proof}
Пусть призмоид $P$ задан системой линейных неравенств $b_i\goe(a_i,x)\goe c_i$, (см. Предложение).

Среди линейных функций $\{a_i\}$ выберем $d$ линейно независимых. Не умаляя общности, пусть это $a_1, \ldots, a_d$. 
Многогранник, задаваемый системой неравенств
$$b_i\goe(a_i,x)\goe c_i, \qquad i = 1, 2, \ldots, d,$$
есть параллелепипед , который обозначим через $\Pi$.  

По Предложению~\ref{prop:1}, для произвольной вершины $v_j$ многогранника $P$ и любого $i = 1, 2, \ldots, d$ выполняется или $(a_i, v)= b_i$ или 
$(a_i, v)= c_i$. Следовательно $v$ -- вершина параллелепипеда $\Pi$. Таким образом, $Vert(P)\subseteq Vert(\Pi)$.

Сделаем аффинное преобразование, переводящее $\Pi$ в единичный куб $\{0,1\}^d$. Тогда каждая вершина многогранника $P$ перейдет в вершину куба. Следовательно,
$P$ перейдет в $0/1$-многогранник.

\end{proof}

\begin{thm}\label{del}
Любой совершенный призмоид аффинно эквивалентен некоторому решетчатому многограннику Делоне.
\end{thm}

\begin{proof}
Пусть $P$ --- $d$-мерный совершенный призмоид. По Теореме~\ref{01}, $P$ аффинно эквивалентен 0/1-многограннику, полученному из
$d$-мерного куба. Это значит, что в $\mathbb R^d$ можно ввести такую аффинную систему координат, в которой все вершины многогранника $P$ --- целые точки.

Пусть $v_1, \ldots, v_k$ --- все вершины $P$. Рассмотрим множество
$$\Lambda(P) = \left \{u : u = \sum\limits_{j = 1}^k  n_j v_j, \;\text{где } n_j \in \mathbb Z,\; \sum\limits_{j = 1}^k n_j = 1 \right \}.$$

Множество $\Lambda(P)$ дискретно, поскольку все его точки целые в выбранной системе координат. Кроме того, параллельные переносы вида
$$m_1 v_1 + \ldots + m_k v_k, \qquad m_i \in \mathbb Z,\; \sum\limits_{i = 1}^k m_i = 0$$
сохраняют $\Lambda(P)$ и образуют группу, транзитивно действующую  на множество $\Lambda(P)$. Поэтому $\Lambda(P)$ --- решетка.

Пусть призмоид $P$ задан системой линейных неравенств $b_i\goe(a_i,x)\goe c_i$ (см. Предложение).

Пусть $q_i(x)=(-b_i+(a_i,x))(-c_i+(a_i,x))$ --- неоднородная функция второго порядка с квадратичной формой ранга 1. Для произвольной точки $u\in\Lambda(P)$ докажем, что 
$q_i(u)\goe 0$.

Действительно: $u = n_1 v_1 + \ldots + n_k v_k$, где  $n_1 + \ldots + n_k = 1$ и $n_j \in \Z$. Тогда:
$$q_i(u)=\left(-b_i + \sum\limits_{j=1}^k n_j(a_i,v_j)\right)\left(-c_i + \sum\limits_{j=1}^k n_j(a_i,v_j)\right).$$
Далее:
$$ \sum\limits_{j=1}^k n_j(a_i,v_j)=(n_{j_1}+\ldots+n_{j_l})b_i+(n_{j_{l+1}}+\ldots+n_{j_k})c_i=pb_i+(1-p)c_i,\text{ где }p=n_{j_1}+\ldots+n_{j_l}\in\Z$$
Пусть $p < 0$. Тогда перепишем: $p b_i + (1 - p) c_i = c_i + (-p) \cdot (c_i - b_i)$. Поскольку $c_i < b_i$, имеем $(-p) \cdot (c_i - b_i) < 0$.
В результате $p b_i + (1 - p) c_i < c_i$ и тем более $p b_i + (1 - p) c_i < b_i$. Следовательно, обе скобки в $q_i(u)$ отрицательны. Таким образом $q_i(u)>0$.

Пусть $p = 0$. Тогда $p b_i + (1 - p) c_i = c_i$. Значит, $q_i(u) = 0$.

Пусть $p = 1$. Тогда $p b_i + (1 - p) c_i = b_i$. Значит, $q_i(u) = 0$.

Пусть $p > 1$. Тогда перепишем: $p b_i + (1 - p) c_i = b_i + (1 - p)(b_i - c_i)$, что верно. Поскольку $c_i < b_i$, имеем $(1 - p) (b_i - c_i) > 0$. 
В результате $p b_i + (1 - p) c_i > b_i$ и тем более $p b_i + (1 - p) c_i < c_i$. Следовательно, обе скобки в $q_i(u)$ положительны. Таким образом $q_i(u)>0$.

Рассмотрим квадратичную функцию
$$Q(x)=\sum\limits_{i=1}^N q_i(x).$$ 
Легко видеть, что $Q(x)$ --- неоднородная функция второго порядка с положительно определенной однородной составляющей.

Для любой точки $u$ решетки $\Lambda(P)$ имеем $q_i(u) \geq 0$, а следовательно, и $Q(u)\goe 0$. 
Более того, если $Q(u) = 0$, то $q_i(u) = 0$ для любого $i = 1, 2, \ldots, N$.

Условие $q_i(u) = 0$ эквивалентно альтернативе: либо $(a_i,u)=b_i$, либо $(a_i,u)=c_i$. Поэтому если $u$ --- вершина многоранника $P$, то, в силу
Предложения, $Q(u) = 0$. Поэтому уравнение $Q(x) = 0$ определяет пустой эллипсоид, описанный вокруг $P$.

Осталось показать, что $Q(u) \neq 0$ при всех $u \in \Lambda(P)\backslash Vert(P)$. Предположим противное: $Q(u) = 0$, $u \in \Lambda(P)\backslash Vert(P)$. Поскольку $u$ лежит на поверхности эллипсоида, описанного вокруг $P$, $u$ находится строго вне $P$. Тогда для некоторого $i$ верно
$(a_i, u) > b_i$. Но тогда $q_i(u) > 0$ и $Q(u) > 0$, противоречие.

Переводя аффинным преобразованием эллипсоид $Q(u) = 0$ в сферу, убеждаемся, что $P$ аффинно эквивалентен решетчатому многограннику Делоне. 

\end{proof}

\end{document}